\documentclass[a4paper,14]{article}

\usepackage{amsmath}
\usepackage{mathtools}
\usepackage{titling}
\usepackage{graphicx}
\usepackage{tikz}
\usepackage{enumitem}
\usepackage[sorting=nty]{biblatex}
\usepackage{caption}
\usepackage[T1]{fontenc}
\usepackage[latin1]{inputenc}
\usepackage{amsthm}
\usepackage{tabto}
\usepackage{float}
\usepackage{amsthm}
\usepackage{thmtools}
\usepackage{xpatch}

\xpatchbibdriver{article}
  {\usebibmacro{title}%
   \newunit}
  {\usebibmacro{title}%
   \printunit{\addcomma\space}}
  {}
  {}
\xpatchbibdriver{book}{%
\newunit\newblock
\printfield{note}%
}
{%
\setunit{\addcomma\space}\newblock
\printfield{note}%
}{}{}

\usetikzlibrary{arrows}
\predate{}
\postdate{}
\newtheorem{theorem}{Theorem}
\newtheorem{lemma}{Lemma}
\definecolor{white}{rgb}{1,1,1}
\declaretheoremstyle[ headfont=\normalfont,bodyfont=\itshape]{ackTitle}
\declaretheorem[shaded={margin=10pt,bgcolor=white}, style=ackTitle, title=\textsc{Acknowledgment}, numbered=no]{acknowledgment}
\declaretheoremstyle[bodyfont=\normalfont, spaceabove=5pt]{normalbody}
\declaretheorem[style=normalbody]{case}
\renewcommand{\mid}{$ such that $}
\renewcommand{\:}{\,}
\begin{document}

\title{Paths in tournaments, a simple proof of Rosenfeld's Conjecture}
\author{Charbel Bou Hanna\\\\
Lebanese University, KALMA Laboratory, Baalbek\\
University of Angers, LAREMA Laboratory, Angers}

\date{}
\maketitle

\begin{abstract}
\noindent
Rosenfeld Conjectured \parencite{rosenfeld} in 1972 that there exists an integer $K\geq 8$ such that any tournament of order $n\geq K$ contains any Hamiltonian oriented path. In 2000, Havet and Thomass\'{e} \parencite{havetthomasse} proved this conjecture for any tournament with exactly 3 exceptions. We give a simplified proof of this fact.
\end{abstract}
\section{Introduction}
A tournament $T$ is an orientation of a complete graph. The set of vertices of $T$ is denoted by $V(T)$ and the set of arcs by $E(T)$. $v(T)$ will denote the order of $T$, $|V(T)|$. Sometimes, we write $|T|$ instead of $v(T)$. $T$ is said to be an $n$-tournament if $v(T)=n$.\\
The out-neighbour (resp. in-neighbour) of a vertex $v$ in $T$ is denoted by $N^+_T(v)$ (resp. $N^-_T(v)$) and its out-degree (resp. in-degree) is denoted by $d^+_T(v)$ (resp. $d^-_T(v)$). We denote by $\delta^+(T)$ (resp. $\delta^-(T)$) the minimal out-degree (resp. in-degree) and by $\Delta^+(T)$ (resp. $\Delta^-(T)$) the maximal out-degree (resp. in-degree). Note that $\delta^-(T)\leq \Delta^+(T)$ (resp. $\delta^+(T)\leq\Delta^-(T)$). A tournament $T$ is said to be regular if $d^+(v)=d^-(v)\,\forall\, v\in V(T)$. A cyclic triangle is a circuit of length 3. We denote by $T^+_4$ a tournament composed of a circuit triangle together with a source. A Paley tournament on 7 vertices is a tournament $T$ such that $V(T)=\{v_i,1\leq i\leq 7\}$ and $(v_i,v_j)\in E(T)$ if and only if $j-i\equiv 1,\,2$ or $4(mod\,7)$. We write $T'\subseteq T$ whenever $T'$ is a subtournament of $T$. Let $S\subseteq V(T)$, we denote by $T[S]$ the subtournament of $T$ induced by $S$. If $T'\subseteq T$ and $S\subseteq V(T)$, we write $T'+S=T[V(T')\cup S]$ and $T'-S=T[V(T')-S]$. Let $v\in V(T)$, $d^+_S(v)=|N^+_T(v)\cap S|$ and $d^-_S(v)=|N^-_T(v)\cap S|$. A subset $\{x_1,...,x_r\}$ in $T$ will be denoted by $[x_1,x_r]$.\\
Let $P=x_1...x_s$ be an oriented path, set $\widetilde{P}=x_s...x_1$. $P$ is called an $s$-path, $x_1$ and $x_s$ are its extremities, $x_1$ is the origin and $x_s$ is the end. The length of $P,\,l(P)$, is the number of its arcs. $P$ is said to be directed if all of its arcs are oriented in the same direction. A block of $P$ is a maximal (for $\subseteq$) directed subpath of $P$. The path $P$ is said  to be of type $P(b_1,...,b_m)$ and we write $P=P(b_1,...,b_m)$, if $P$ is composed of $m$ successive blocks $B_1,...,B_m$, such that $l(B_i)=b_i$. Moreover, we write $P=P^+(b_1,...,b_m)$, if $(x_1,x_2)\in E(P)$. Else, we write $P=P^-(b_1,...,b_m)$. $P$ is said to be antidirected if each block of $P$ is of length 1.\\
Note that if $P=x_1...x_s=P^+(b_1,...,b_m)$ (resp. $P^-(b_1,...,b_m)$) and $P'=x'_1...x'_s=P^+(b_1,...,b_m)$ (resp. $P^-(b_1,...,b_m)$), then $P$ and $P'$ are isomorphic. We write $P\equiv P'$. Furthermore, if we write $x_1...x_s\equiv x'_1...x'_s$, then the mapping: 
\begin{center}
\begin{tabular}{l c c c c c}
$f:$&$V(P)$&$\rightarrow$&$V(P')$  & is an isomorphism.\\
&$x_i$&$\rightarrow $&\multicolumn{1}{l}{$f(x_i)=x'_i $}&
\end{tabular}\\
\end{center}
A path $P$, in a tournament $T$, is said to be Hamiltonian if $V(P)=V(T)$. Let $P=x_1...x_s$ and $Q=y_1...y_r$ be 2 disjoint paths in $T$, $PQ$ denotes the path $x_1...x_s\, y_1...y_r$. In a similar way, we may define the path $P_1P_2...P_t$ from $t$ pairewisely disjoint paths $P_1,...,P_t$ of $T$. Let $P=v_1...v_{i-1}\, v_i\, v_{i+1}...v_n$ be an oriented path in a tournament $T$, then $P-v_i$ is the path $v_1..v_{i-1}\, v_{i+1}...v_n$ in $T$.\\
A strong tournament is such that any two of its vertices can be joined by a directed path. It is known that any tournament $T$ is a transitive union of strong subtournaments (called strong components) $I_1,\,I_2,...,I_t$ that is $(v_i,v_j)\in E(T)\,\forall\, i<j$, $v_i\in I_i$ and $v_j\in I_j$. We write $T=I_1...I_t$. A strong tournament is characterized by the following property due to Camion \parencite{camion}.
\begin{theorem}
$T$ is strong if and only if $T$ contains a Hamiltonian circuit.
\end{theorem}
\noindent
As a consequence of the above theorem, if $T=I_1...I_t$, then $T$ contains a directed path ending at $x$ with $V(P)=V(I_1)\cup...\cup V(I_i)$ for every $x\in T_i$.\\
The complement of a digraph $D$, denoted by $\overline{D}$, is the digraph obtained from $D$ by reversing all its arcs.\\
\\
In 1971, Gr\"{u}nbaum \parencite{grunbaum} proved the following theorem:
\begin{theorem} Any tournament contains any Hamiltonian antidirected path with exactly 3 exceptions: a cyclic triangle ($T_3$), a regular tournament on 5 vertices ($T_5$) and a Paley tournament on 7 vertices ($T_7$).
\end{theorem}
\noindent
Set $\mathcal{T}_{3,5,7}=\{T_3,T_5,T_7\}$. Note that if $T \in \mathcal{T}_{3,5,7}$ contains a copy of a path $P$, then any vertex in $T$ is an origin of a copy of $P$. Rosenfeld \parencite{rosenfeld}, in 1972, inspired by the work of Gr\"{u}nbaum, conjectured that there exists $K\geq 8$ such that any tournament of order $n\geq K$ contains any Hamiltonian oriented path. The case of directed path being R\'{e}dei's theorem \parencite{redei}. Alspach, Rosenfeld \parencite{alspach} and Straight \parencite{straight} proved Rosenfeld's conjecture on paths of 2 blocks. In 1973, Forcade \parencite{forcade} proved Rosenfeld's conjecture for any tournament of order $2^n$. Thomason \parencite{thomason} was the first one to give a general answer. He proved, in 1986, that there exists $n_0<2^{128}$ such that for all $n\geq n_0$ any tournament of order $n$ contains any Hamiltonian oriented path. Havet and Thomass\'{e} settled the problem by proving that the three exceptions of Gr\"{u}nbaum are the only tournaments not satisfying Rosenfled's conjecture. Havet and Thomass\'{e}'s proof consists of giving a refinement of a key idea introduced by Thomason saying that any set of $b_1+1$ vertices in an $n$-tournament contains an origin of any $(n-1)$-path whose first block is of length $b_1$. They proved that if $s^+(x,y)\geq b_1+1$, then $x$ or $y$ is an origin of a copy of such path where $s^+(x,y)=|\{z\in T$ such that $ z$ can be reached from $x$ or $y$ by a directed path$\}|$. This new performance allowed them to remark that proving the existence of an $(n-1)$-oriented path in any $n$-tournament $T$ is equivalent to the existence of any Hamiltonian path $P$ in this tournament unless $(T,P)$ is one of the 69 exceptions that were verified one by one. In this paper, we give a simple proof of the result without treating all these exceptions.
\section{The main result}
In our proof, we are going to use the following theorem due to El-sahili and Ghazo-Hanna \parencite{sahilighazo}.
\begin{theorem}
A path lies in a tournament $T$ if and only if it lies in $\overline{T}$.\\
\\{\normalfont The following Lemmas will be useful in the sequel.}
\begin{lemma}
Let $T$ be a tournament then $|\{ v\in T\,/\, T-v\in \mathcal{T}_{3,5,7}\}|\leq 2$.
\begin{proof}
Suppose that $T$ contains 3 distinct vertices, $v_1$, $v_2$, $v_3\mid T-v_i\in \mathcal{T}_{3,5,7}$, $1\leq i\leq 3$. We may suppose without loss of generality that $(v_1,v_2)$, $(v_2,v_3)\in E(T)$. Since $T-v_1\in \mathcal{T}_{3,5,7}$ and $T-v_3\in \mathcal{T}_{3,5,7}$, then $N^-_{T-v_3}(v_1)=N^-_{T-v_1}(v_3)$ and $N^+_{T-v_1}(v_3)=N^+_{T-v_3}(v_1)$, but $v_2\in N^+_{T-v_3}(v_1)$ and $v_2\notin N^+_{T-v_1}(v_3)$, which is a contradiction.
\end{proof}
\noindent
{\normalfont We may deduce from this Lemma more precise conclusions. Indeed, if $T\in \mathcal{T}_{3,5,7}$, then $\forall\:x\in T$, $\exists\, y\in N^+(x)\mid T-\{x,y\}\notin \mathcal{T}_{3,5,7}$}.
\end{lemma}
\begin{lemma} (Simple Lemma) Let $n\geq 4$. Suppose that any $s$-tournament $(s<n)$ contains any $s$-path unless if the tournament is in $\mathcal{T}_{3,5,7}$ and the path is antidirected. Let $T$ be an $n$-tournament with a vertex $v\mid d^-(v)=0$ and let $P=x_1...x_n$ be a non directed path. Unless $P=P^+(1,2)$ and $T-v$ is a cyclic triangle, we have:
\begin{enumerate}[label=\textnormal{\arabic*.}, leftmargin=\parindent, labelindent=\parindent, listparindent=\parindent,]
\item $T$ contains a copy of $P$ with origin $x\neq v$.
\item Any of the vertices of $T-v$ is an origin of a copy of $P$ if $T-v\in \mathcal{T}_{3,5,7}$.
\end{enumerate}
\begin{proof}
\quad \\
\vspace*{-\baselineskip}
\begin{enumerate}[label=\textnormal{\arabic*.}]
\item Let $j\in [2,n]$ be the minimal integer $\mid d^-_P(x_j)=0$. If $j=n$, then $T-v$ contains a path $P'$ such that $P'v\equiv P$ since otherwise $P-x_j$ is antidirected and so $P=P^+(1,2)$ and $T=T^+_4$ which is a contradiction. Otherwise, $1<j<n$. In this case, $P-x_j+(x_{j-1},x_{j+1})$ or $P-x_j+(x_{j+1},x_{j-1})$ is not antidirected. Suppose that $P'=P-x_j+(x_{j-1},x_{j+1})$ is not antidirected, then $T-v$ contains a path $v_1\,v_2...v_{j-1}\,v_{j+1}...v_n\equiv P'$ and so $v_1...v_{j-1}\:v\:v_{j+1}...v_n\equiv P$ and $v_1\neq v$.
\item If $T-v\in \mathcal{T}_{3,5,7}$, then we may suppose that a copy of $P'$ may be found in $T-v$ starting at any one of the vertices of $T-v$.
\end{enumerate}
\end{proof}
\noindent
{\normalfont By analogy, the above lemma is valid if $d^+(v)=0$ with the only exception where $P=P^-(1,2)$ and $T-v$ is a cyclic triangle. Using the same reasoning, we may deduce that if $P$ contains a vertex $x_j$ with $1<j<n\mid d^-_P(x_j)=0$, then $T$ contains a copy of $P$ such that $v$ is not an extremity of $P$.\\
We are going now to present the proof of the main result.}
\end{lemma}
\end{theorem}
\begin{theorem}
Any $n$-tournament contains any $n$-path unless the three exceptions of Gr\"{u}nbaum.\\
\begin{proof}
Let $T$ be an $n$-tournament and let $P=x_1..x_n$ be an $n$-path. If $n\leq 4$, then any $n$-path is either directed, antidirected or of two blocks. By theorem 3, the problem will be solved if we prove the existence of \^{P} in \^{T}, where \^{P}$\in \{P,\overline{P}\}$ and \^{T}$\in\{T,\overline{T}\}$. We may suppose, without loss of generality, that $\Delta^+(T)\geq \Delta^+(\overline{T})$. We argue by induction on $|T|\geq 5$. The case $\delta^-(T)=0$ can be deduced from simple lemma and R\'{e}dei's theorem. Set $\delta^-(T)=i\geq 1$ and suppose, without loss of generality, that $(x_i,x_{i+1})\in E(P)$, since otherwise, we use $\overline{P}$. Then $x_i...x_n$ is of type $P^+(b_1,b_2,...,b_r)$ for some $b_1,...,b_r\geq 1$. Let $v\in T\mid d^-(v)=\delta^-(T)$, $T_1=T[N^-(v)]$ and $T_2=T[N^+(v)]$. We will treat the two cases according to the existence of a copy of $x_1...x_i$ in $T_1$.
\begin{case}
$T_1\not\supseteq v_1...v_i\equiv  x_1...x_i$.\\
\noindent
In this case, $x_1...x_i$ is antidirected and $T_1\in \mathcal{T}_{3,5,7}$. If $b_1\geq 3$, let $a\in T_2\mid T_2-a\supseteq v_{i+3}...v_n \equiv x_{i+3}...x_n$. Then $T_1+a\supseteq v_1...v_{i+1}\equiv x_1...x_{i+1}\mid v_{i+1}\neq a$. So $v_1...v_{i+1}\:v\:v_{i+3}...v_n\equiv P$.\\
In the case $b_1\leq 2$, \textbf{for} $\mathbf{(x_i,x_{i-1})\in E(P)}$, suppose that $b_1=1$, we discuss if $T_2\supseteq v_{i+2}...v_n\equiv x_{i+2}...x_n$ or not. In the first case, suppose that there exists $x\in T_1$ such that $(v_{i+2},x)\in E(T)$. By simple lemma, $T_1+v\supseteq v_1...v_{i+1}\equiv x_1...x_{i+1}$ with $v_{i+1}=x$, so $v_1...v_n\equiv P$. If $(x,v_{i+2})\in E(T)\,\forall\, x\in T_1$, let $t$ be the minimal integer $(t\geq i+2)$ such that $d^-_P(x_t)=0$. If there exists $y\in T_1$ such that $T_2+y\supseteq w_{i+1}...w_n\equiv x_{i+1}...x_n$ with $w_{i+1}=v_{i+2}$, then $(T_1-y)+v\supseteq v_1...v_i\equiv x_1...x_i$ and so $v_1...v_i\,w_{i+1}...w_n\equiv P$. Otherwise, $N^-(x)\cap[v_{i+2},v_t]=\phi$ and $(v_{t+1},x)\in E(T)\,\forall\, x\in T_1$. In this case, let $y\in T_1$. $Q=(v_{t+1}\,y\,v_t\,v_{t-1}...v_{i+2}-v_t)-v_{i+2}$ is a directed path. If $t\geq i+3$, then $l(Q)=b_2-1$. $T_2[v_{t+2},v_n]+\{v_t,v\}\supseteq w_t...w_n\equiv x_t...x_n$ with $w_t=v$ or, by simple lemma, $w_t=v_t$. In the other hand, by simple lemma, $T_1+v_{i+2}\supseteq v_1...v_i\, v_{i+1}\equiv x_1...x_{i+1}$ with $v_{i+1}=y$ and $v_i\neq v_{i+2}$, so $v_1...v_i\, \widetilde{Q}\, w_t...w_n\equiv P$. If $t=i+2$, let $x\in T_1, T_1-x\supseteq v_2...v_i\equiv x_2...x_i$. If $(v_{t+1},v_{t+2})\in E(T), v_{t+1}\,v_2...v_i\, v_t\,x\,v\, v_{t+2}...v_n\equiv P$. If $(v_{t+2},v_{t+1})\in E(T)$ and there exists $y\in T_1$ such that $(v_{t+2},y)\in E(T)$, consider $(a,b)\in E(T_1-y)$, $(T_1-\{a,b,y\})+v\supseteq v_1...v_{i-2}\equiv x_1...x_{i-2}$ and so $v_1...v_{i-2}\, v_t\,a\,b\,v_{t+1}\,y\,v_{t+2}..v_n\equiv P$. Else, let $z\in T_1$. $(T_1-z)+v_{t+1}\supseteq v_1...v_i\equiv x_1...x_i$ with $v_i\neq v_{t+1}$ and $(T[v_{t+2},v_n])+v\supseteq w_{t+1}...w_n\equiv x_{t+1}...x_n$ with $w_{t+1}=v_{t+2}$ then $v_1...v_i\,v_t\,z\,w_{t+1}...w_n\equiv P$.\\
In the other case, $T_2\not\supseteq v_{i+2}...v_n\equiv x_{i+2}...x_n$, so $T_2\in \mathcal{T}_{3,5,7}$ and $x_{i+2}...x_n$ is antidirected. If $b_2=1$, then $P$ is antidirected treated in Theorem 2. Otherwise, $b_2=2$. Let $a\in T_2$, then $T_1+a\supseteq v_1...v_{i+1}\equiv x_1...x_{i+1}$ such that $v_{i+1}\neq a$. Since $v_{i+1}\in T_1\in \mathcal{T}_{3,5,7}$, then $d^-_{T_2}(v_{i+1})\geq 2$, there exists $b\in T_2-a$ such that $(b,v_{i+1})\in E(T)$. If $T_2-\{a,b\}\supseteq Q\equiv x_{i+4}...x_n$, then $v_1...v_{i+1}\,b\,v\,Q\equiv P$. Otherwise, there exists $c\in T_2-\{a,b\}$ such that $(c,b)\in E(T)$, then, by simple lemma, $(T_2-\{a,b\})+v\supseteq Q\equiv x_{i+3}...x_n$ starting by c, then $v_1...v_{i+1}\,b\,Q\equiv P$.\\
\textbf{For} $\mathbf{b_1=2}$, here also we study if $T_2\supseteq v_{i+2}...v_n\equiv x_{i+2}...x_n$ or not. In the first case, let $x\in N^-(v_{i+2})\cap T_1$, if any, by simple lemma, $T_{1}+v\supseteq v_1...v_{i+1}\equiv x_1...x_{i+1}$ with $v_{i+1}= x$ and so $v_1...v_n\equiv x_1...x_n$. If $N^-(v_{i+2})\cap T_1=\phi$, $v$ can be inserted inside $v_{i+3}...v_n$ to obtain a path $Q\equiv x_{i+2}...x_n$ starting at $v_{i+3}$ in $(T_2-v_{i+2})+v$. Similarly, by simple lemma, if $N^-(v_{i+3})\cap T_1=\phi$, we may find a copy of $P$ in $T$, so, we may suppose in the sequel that $d^-_{T_1}(v_{i+2})=d^-_{T_1}(v_{i+3})=0$. As above, let $t\geq i+3$ be the minimal integer such that $d^-_P(x_t)=0$. If $t=i+3$, since $T_1\in \mathcal{T}_{3,5,7}$, then $d^+_{T_2}(x)\geq 1\,\forall\,x\in T_1$, since otherwise, $\Delta^+(\overline{T})>\Delta^+(T)$. If there exists $x\in T_1$ and $a\in T_2$ such that $(x,a)\in E(T)$ and $d^-_{T_2-v_{i+2}}(a)\geq 1$, let $b\in N^-_{T_2-v_{i+2}}(a)$. By simple lemma, $T_1+v_{i+2}\supseteq v_1...v_{i+1}\equiv x_1...x_{i+1}$ with $v_{i+1}=x$. $(T_2-\{a,v_{i+2}\})+v\supseteq v_{i+3}...v_n\equiv x_{i+3}...x_n$ with $v_{i+3}=v$ or, by simple lemma, $v_{i+3}=b$, so $v_1...v_{i+1}\, a\,v_{i+3}...v_n\equiv P$. Otherwise, if $(v_{i+2},a)\in E(T)$, then the problem is solved by considering $v_1...v_{i+1}\equiv x_1...x_{i+1}$ in $T_1+v_{i+3}$. Otherwise, there exists $a\in T_2$ such that $\forall\, x\in T_1$, $\forall\,y\in T_2-a$, we have $(x,a)$, $(y,x)$, $(a,y)\in E(T)$. Let $y\in T_2-a$ such that $T_2-\{a,y\}\notin T_{3,5,7}$. $T_1+y\supseteq v_1...v_{i+1}\equiv x_1...x_{i+1}$ with $v_{i+1}\neq y$ and $T_2-\{a,y\}\supseteq v_{i+4}...v_n\equiv x_{i+4}...x_n$ then $v_1...v_{i+1}\,a\,v\,v_{i+4}...v_n\equiv P$. If $t>i+3$, let $x\in T_1$, $T_1-x\supseteq v_2...v_i\equiv x_2...x_i$. If $(T_2-v_{i+2})+\{v,x,v_i\}\supseteq w_i...w_n\equiv x_i...x_n$ with $w_i\in \{v_{i+3},v_i\}$, then $v_{i+2}\, v_2...v_{i-1}\,w_i...w_n\equiv P$. Otherwise, $\forall\, x\in T_1$ $N^-(x)\cap [v_{i+4},v_t]=\phi$ and $(v_{t+1},x)\in E(T)$. Let $x,\,y\in T_1$ such that $(x,y)\in E(T)$ and $T_1-\{x,y\}\supseteq v_2...v_{i-1}\equiv x_2...x_{i-1}$. As above, $(T_2-\{v_{i+2},v_{i+3}\})+\{y,v\}\supseteq w_{i+2}...w_n\equiv x_{i+2}...x_n$ with $w_{i+2}\in \{y,v_{i+4}\}$, so $v_{i+2}\,v_2...v_{i-1}\,v_{i+3}\,x\, w_{i+2}...w_n\equiv P$. If $T_2\not\supseteq v_{i+2}...v_n\equiv x_{i+2}...x_n$ then, $x_{i+2}...x_n$ is antidirected and $T_2\in \mathcal{T}_{3,5,7}$. If $|T_2|>|T_1|$, the problem is solved by considering $\overline{\widetilde{P}}\, (b_1=1)$. Otherwise, $|T_1|=|T_2|$. Let $x\in T_1$. $T_2+x \supseteq v_{i+1}...v_n\equiv x_{i+1}...x_n\mid v_{i+1}\neq x$. If there exists $y\in T_1-x$ such that $(y,v_{i+1})\in E(T)$, then $(T_1+v)-x \supseteq w_1..w_i\equiv x_1...x_i$ with $w_i=y$ and so $w_1...w_i\:v_{i+1}..v_n\equiv P$, unless $T_1-\{x,y\}\in \mathcal{T}_{3,5,7}$, then there exists $z\in T_1-\{x,y\}$ such that $(y,z)\in E(T)$. By simple lemma, $(T_1-\{x,y\})+v\supseteq w_1...w_{i-1}$ with $w_{i-1}=z$, and so $w_1...w_{i-1}\, y\, v_{i+1}...v_n\equiv P$. Otherwise $T_1$ and $T_2$ are cyclic triangles. The problem is solved unless if for $u\in T,\:T[N^+(u)]$, $T[N^-(u)]\in T_3$. So, if $u\in T_1$, then $d^+_{T_2}(u)=1$. Else, $d^+_{T_1}(u)=2$. Set $V(T_1)=\{x,y,z\}$, then $T_2+x\supseteq v_4...v_7\equiv x_4...x_7$. Since $x_4...x_7$ is antidirected, then $v_5=x$. Suppose, without loss of generality that, $(y,v_4)\in E(T)$, then $z\:v\:y\:v_4...v_7\equiv P$.\\
\textbf{If} $\mathbf{(x_{i-1},x_i)\in E(P)}$, let $x\in T_1$. $T_1-x \supseteq Q_1\equiv x_1...x_{i-1}$. If $T_2+x\supseteq v_{i+1}...v_n\equiv x_{i+1}...x_n$ and $v_{i+1}\neq x$, $Q_1\:v\: v_{i+1}...v_n\equiv P$. If $v_{i+1}=x$, then, by simple lemma, $T_1+v \supseteq v_1...v_{i+1}\equiv x_1...x_{i+1}\mid v_{i+1}=x$. Then $v_1...v_n\equiv P$. If $T_2+x$ contains no copy of $x_{i+1}...x_n$, then this path is antdirected, so we consider $\widetilde{P}=y_1...y_n$ to remark that $P$ or $\overline{P}$ contains the arcs $(y_i,y_{i-1})$ and $(y_i,y_{i+1})$, we recover a previous case.
\end{case}
\begin{case}
$T_1\supseteq v_1...v_i\equiv x_1...x_i$.\\
If $b_1\geq 2$, then either $T_2\supseteq v_{i+2}...v_n\equiv x_{i+2}...x_n$ or not. In the first case, $v_1...v_i\:v\:v_{i+2}...v_n\equiv P$ and in the last case, $T_2\in \mathcal{T}_{3,5,7}$ and $x_{i+2}...x_n$ is antidirected. If $\exists\: x\in T_2\mid (v_i,x)\in E(T)$, then, by simple lemma, $T_2+v\supseteq v_{i+1}...v_n\equiv x_{i+1}...x_n$ with $v_{i+1}=x$. Thus $v_1...v_i\:v_{i+1}...v_n\equiv P$. Otherwise, $(x,v_i)\in E(T)\:\forall\: x\in T_2$, we have $(v_i,v_j)\in E(T)$ whenever $j<i$, since otherwise $\Delta^+(\overline{T})\geq d^+_{\overline{T}}(v_i)=d^-_T(v_i)>|T_2|=d^+(v)=\Delta^+(T)$, a contradiction. By simple lemma, $T_2+v_i$ contains a path $v_{i+1}...v_n\equiv x_{i+1}...x_n$ with $v_{i+1}\in T_2$. If $i=1$, then $v\,v_{i+1}...v_n\equiv P$. Otherwise, let $j<i$ be the maximal integer such that $d^+_P(x_j)=0$, then $v_1...v_{j-1}\,v\,v_j...v_{i-1}\equiv x_1...x_i$. Since $(v_i,v_{i-1}) \in E(T)$, then $\exists\,x\in T_2\mid  (v_{i-1},x)\in E(T)$. By simple lemma, we suppose that $v_{i+1}=x$ and so $v_1...v_{j-1}\,v\,v_j...v_{i-1}$ $v_{i+1}...v_n\equiv P$.\\
Now, we will study the case $b_1=1$. Suppose that $N^+(v_i)\cap T_2\neq\phi$. Let $I_1,\,I_2,...I_t$ be the strong connected components of $T_2$ such that $T_2=I_1...I_t$ and let $l=max\{j, N^+(v_i)\cap I_j\neq \phi\}$, $I=I_1\,I_2...I_l$ and $s=|I|$. We discuss according to the value of $s$:\\
For the case $\mathbf{s>b_2}$, let $Q=u_1...u_{b_2+1}$ be a directed path in $I\mid u_{b_2+1}\in N^+(v_i)$ and let $j=i+b_2+2$. So either $T_2-[u_2,u_{b_2+1}]$ contains a path $v_j...v_n\equiv x_j...x_n$, and in this case $v_1...v_i\:u_{b_2+1}...u_2\:v\:v_j...v_n\equiv P$, or $T_2-[u_2,u_{b_2+1}]\in \mathcal{T}_{3,5,7}$ and $x_j...x_n$ is antidirected. By simple lemma, $(T_2-[u_2,u_{b_2+1}])+v$ contains a path $v_{j-1}...v_n\equiv x_{j-1}...x_n$ with $v_{j-1}=u_1$. In this case, $v_1...v_i\,u_{b_2+1}...$ $u_2\,v_{j-1}...v_n\equiv P$.\\
Now we will study the case $\mathbf{s=b_2}$. Let $Q=u_1...u_{b_2}$ be a directed path in $I\mid u_{b_2}\in N^+(v_i)$ and let $j=i+b_2+2$. If $T_2-I\notin \mathcal{T}_{3,5,7}$ or $x_j...x_n$ is not antidirected, then $v_1...v_i\,u_{b_2}...u_1\,v\,v_j...v_n\equiv P$ where $T_2-I\supseteq v_j...v_n\equiv x_j...x_n$. Otherwise, we will continue the proof depending on the orientation of $x_i\,x_{i-1}$ and on the value of $i$. If $(x_{i-1},x_i)\in E(P)$ or $i=1$, let $a\in T_2-I$. By simple lemma, $(T_2-(I\cup a))+v_i$ contains a path $v_j...v_n\equiv x_j...x_n$ with $v_j\neq v_i$, then $v_1...v_{i-1}\: v\: a\: u_{b_2}...u_1\: v_j...v_n\equiv P$. Otherwise, $(x_i,x_{i-1})\in E(P)$ and $i\geq 2$. If $N^+(v_{i-1})\cap (T_2-I)\neq\phi$, let $a\in N^+(v_{i-1})\cap (T_2-I)$. As above $\exists\:r<i\mid v_1...v_{r-1}\:v\:v_r...v_{i-1}\equiv x_1...x_i$ and by simple lemma $(T_2-(I\cup a))+v_i\supseteq v_j...v_n\equiv x_j...x_n$ with $v_j\neq v_i$, so $v_1...v_{r-1}\:v\:v_r...v_{i-1}\:a\:u_{b_2}...$ $u_1\:v_j...v_n\equiv P$. Otherwise, let $abc$ be a cyclic triangle in $T_2-I$. Then, if $b_2\geq 2$, $v_1...v_{i-1}\:a\:v_i\:b\:u_{b_2}...u_3\: v\: v_j...v_n\equiv P$ with $v_j...v_n\subseteq (T_2-(I\cup\{a,b\}))+\{u_1,u_2\}\mid v_j...v_n\equiv x_j...x_n$. If $b_2=b_3=1$, $v_1...v_{i-1}\:a\:v_i\: b\: c\:v\:v_{j+2}...v_n\equiv P$ where $v_{j+2}...v_n\subseteq T_2-\{a,b,c\}$. If $b_2=1$, $b_3=2$ and $T_2-I\in \{T_5,T_7\}$, choose $d\in T_2-\{u_1,a,b,c\}$ such that $T_2-\{u_1,a,b,c,d\}\notin \mathcal{T}_{3,5,7}$ ($d$ exists by Lemma 1), then $v_1...v_{i-1}\:d\:v_i\: a\: b\: c\: v$ $v_{j+3}...v_n\equiv P$ where $v_{j+3}...v_n\subseteq T_2-\{a,b,c,d\}$.\\
For the remaining cases, $T_2=T^+_4$. The possible situations are the following. If $|T_1|=2$, then a copy of $P$ is found by considering $\overline{\widetilde{P}}$. If $|T_1|=3$, we have the following two cases:\\
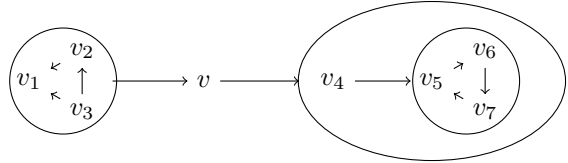
\begin{figure}[H]
\centering
\begin{tikzpicture}[node distance=0.2cm]
\node (v1) at (1,0) {$v_1$};
\node (v2) at (1.7,0.4){$v_2$};
\node (v3) at (1.7,-0.4){$v_3$};
\node (v) at (3.3,0) {$v$};
\node (v4) at (5,0) {$v_4$};
\node (v5) at (6.3,0) {$v_5$};
\node (v6) at (7,0.4){$v_6$};
\node (v7) at (7,-0.4){$v_7$};
\draw[->](v2) to (v1);
\draw[->](v3) to (v1);
\draw[->](v3) to (v2);

\draw[->](v5) to (v6);
\draw[->](v6) to (v7);
\draw[->](v7) to (v5);
\draw[->] (2.1,0) to (v);
\draw[->] (v) to (4.55,0);
\draw[->] (v4) to (6.02,0);
\draw (1.45,0) circle (20pt);
\draw (6.75,0) circle (20pt);
\draw (6.3,0) ellipse (50pt and 30pt);

\end{tikzpicture} 
\captionsetup{labelformat=empty}
\caption{$v_5\:v_4\:v\:v_7\:v_6\:v_3\:v_1\:v_2\equiv P$.}
\end{figure}

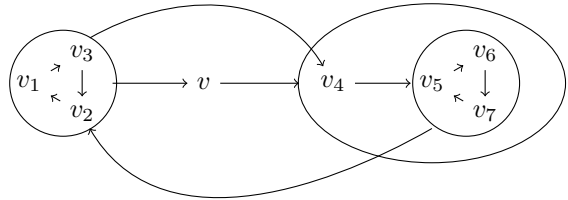
\begin{figure}[H]
\centering
\begin{tikzpicture}[node distance=0.2cm]
\node (v1) at (1,0) {$v_1$};
\node (v3) at (1.7,0.4){$v_3$};
\node (v2) at (1.7,-0.4){$v_2$};
\node (v) at (3.3,0) {$v$};
\node (v4) at (5,0) {$v_4$};
\node (v5) at (6.3,0) {$v_5$};
\node (v6) at (7,0.4){$v_6$};
\node (v7) at (7,-0.4){$v_7$};
\draw[->](v1) to (v3);
\draw[->](v3) to (v2);
\draw[->](v2) to (v1);

\draw[->](v5) to (v6);
\draw[->](v6) to (v7);
\draw[->](v7) to (v5);
\draw[->] (2.1,0) to (v);
\draw[->] (v) to (4.55,0);
\draw[->] (v4) to (6.02,0);
\draw[->] (1.8,0.6) [out=30, in=120] to (v4);
\draw[->] (6.3,-0.6) [out=210, in=300] to (1.8,-0.6);
\draw (1.45,0) circle (20pt);
\draw (6.75,0) circle (20pt);
\draw (6.3,0) ellipse (50pt and 30pt);

\end{tikzpicture}
\captionsetup{labelformat=empty}
\caption{$v_7\:v_6\:v_5\:v_2\:v_3\:v\:v_4\:v_1\equiv P$.}
\end{figure}
\noindent
If $|T_1|=4$, the problem is solved unless $T_1=\overline{T_2}$ and $x_1...x_4\equiv x_9...x_6$.\\
In this case, we have:

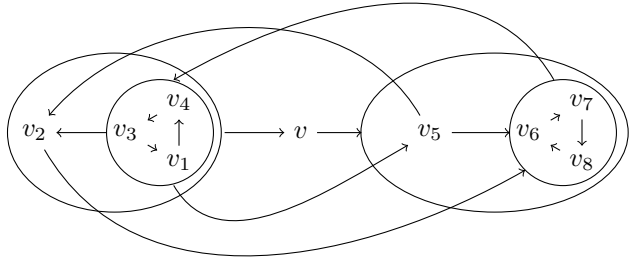
\begin{figure}[H]
\centering
\begin{tikzpicture}[node distance=0.2cm]
\node (v4) at (1.5,0) {$v_2$};
\node (v1) at (2.7,0) {$v_3$};
\node (v3) at (3.4,0.4){$v_4$};
\node (v2) at (3.4,-0.4){$v_1$};
\node (v) at (5,0) {$v$};
\node (v5) at (6.7,0) {$v_5$};
\node (v6) at (8.0,0) {$v_6$};
\node (v7) at (8.7,0.4){$v_7$};
\node (v8) at (8.7,-0.4){$v_8$};
\draw[->](v1) to (v2);
\draw[->](v2) to (v3);
\draw[->](v3) to (v1);

\draw[->](v6) to (v7);
\draw[->](v7) to (v8);
\draw[->](v8) to (v6);

\draw[->] (2.45,0) to (v4);
\draw[->] (3.33,-0.7) [out=300, in=210] to (v5);
\draw[->] (8.33,0.7) [out=120, in=30] to (3.33,0.7);
\draw[->] (4,0) to (v);
\draw[->] (v) to (5.8,0);
\draw[->] (v5) to (7.75,0);
\draw[->] (v4) [out=300, in=210] to (7.95,-0.5);
\draw[->] (v5) [out=120, in=50] to (v4);

\draw (3.15,0) circle (20pt);
\draw (8.45,0) circle (20pt);

\draw (2.55,0) ellipse (40pt and 30pt);
\draw (7.55,0) ellipse (50pt and 30pt);

\end{tikzpicture}
\captionsetup{labelformat=empty}
\caption{$v_3\:v\:v_2\:v_1\:v_4\:v_6\:v_7\:v_8\:v_5\equiv P$.}
\end{figure}
\noindent
Finally, $\mathbf{s\leq b_2-1}$. Let $u_1...u_s$ be a directed path in $T_2 \mid (v_i,u_s)\in E(T)$. If $s\geq 2$, then $b_2\geq 3$. Let $Q$ be a directed path of length $b_2-s-1$ in $T_2-[u_1,u_s]$ such that $(T_2-Q)-[u_1,u_s]\supseteq Q_1\equiv x_{i+b_2+2}...x_n$, and let $Q_2$ be a Hamiltonian directed path in $T[\{u_1,...,u_{s-2},u_s\}]$. If $(x_{i-1},x_i)\in E(P)$ or $|T_1|=1$ then $v_1...v_{i-1}  v  u_s  v_i  Q  u_{s-1}...u_1  Q_1\equiv P$. Otherwise, $(x_i,x_{i-1})\in E(P)$, then $(T_1-v_i)+v \supseteq v'_1...v'_i\equiv x_1...x_i$ such that $v'_i=v_{i-1}$. We will insert $Q,\,v,\,v_i,\, u_s$ into $v_1...v_{i-1}\,u_{s-1}...u_1$ $Q_1$ to obtain a copy of $P$ according to the arcs between $v_{i-1},\,v_i,\,u_{s-1}$ and $u_s$. $(v_{i-1},u_s)\in E(T)$ $\Rightarrow$ $v'_1...v'_i  \,u_s \, v_i \, Q  \,u_{s-1}...u_1\,Q_1\equiv P$. $(u_s,v_{i-1}),$ $(u_{s-1},v_i)\in E(T)$ $\Rightarrow$ $v_1...v_{i-1}\,u_s\,Q\,v\,v_i\,u_{s-1}...u_1\,Q_1\equiv P$. $(u_s,v_{i-1})$, $(v_i,u_{s-1})$,$(u_{s-1},v_{i-1})\in E(T)$ $\Rightarrow$ $v_1...v_{i-1}\,u_{s-1}\,u_s\,v_i\,Q\,u_{s-2}...u_1\,v\,Q_1\equiv P$. The only remainder case is $(u_s,v_{i-1})$, $(v_i,u_{s-1})$, $(v_{i-1},u_{s-1})\in E(T)$. In this case, $v'_1...v'_i\,u_{s-1}\,v_i\,Q\,Q_2\,Q_1\equiv P$. Now, we will treat the case $s=1$ and $b_2\geq 2$. If $(x_i,x_{i-1})\in E(P)$, then $|T_2|\geq 3$. If $\exists\, a\in (T_2-u_1)\cap N^-(v_{i-1})$and $b_2=\Delta^+(T)$, let $Q$ be a Hamiltonian directed path in $T_2-\{u_1,a\}$, then $v_1...v_{i-1}\, a\, v_i\, Q\, u_1\, v\equiv P$. if $b_2<\Delta^+(T)$, let $Q_1$ be a directed path of length $b_2-2$ in $T_2-\{a,u_1\}$ then $T_2-(V(Q_1)\cup\{a\})\supseteq Q_2$ such that $v_1...v_{i-1}\,a\, v_i\, Q_1\, v\, Q_2\equiv P$. If $(T_2-u_1)\cap N^-(v_{i-1})=\phi$, then $(T_1-v_{i-1})\subseteq N^-(v_{i-1})\cap N^-(u_1)$ and $(u_1,v_{i-1})\in E(T)$. Since $b_2\geq 2$, then $x_{i+2}...x_n$ is neither of type $P^+(1,2)$ nor directed starting from $x_{i+2}$. By simple lemma, $(T_2-u_1)+v\supseteq v_{i+2}...v_n\equiv x_{i+2}...x_n$ such that $v_{i+2}\neq v$. For $|T_1|>2$, if $\exists\,j<i-1$ such that $(v_j,v_i)\in E(T)$, then since $(T_1-\{v_{i-1},v_i,v_j\})+u_1\supseteq Q\equiv x_1...x_{i-2}$, so $Q\, v_{i-1}\, v_j\, v_i\,v_{i+2}...v_n\equiv P$. If $(v_i,v_j)\in E(T)$ whenever $j<i$, then $T_1-\{v_1,v_{i-1}\}\supseteq Q_1\equiv x_1...x_{i-2}$ and $Q_1\, v_{i-1}\, v_1\equiv x_{1}...x_i$. Then the problem is solved since $(v_i,v_1)\in E(T)$. For $|T_1|=2$, if $\exists\, a\in T_2-u_1$ such that $(T_2-a)\supseteq v_{i+3}...v_n\equiv  x_{i+3}...x_n$ with $v_{i+3}\neq u_1$, then $a\, v_1 \,v\, v_2\, v_{i+3}\, v_n\equiv P$. Otherwise, by simple lemma, $x_5...x_n$ is a directed path. Let $a,\,b\in T_2-u_1$ such that $(a,b)\in E(T)$. $|T_2|>3$, since otherwise, $N^+_T(b)=\{v_2\}$ then $d^+_T(b)=1$ and so $\Delta^+(\overline{T})>\Delta^+(T)$, a contradiction. Let $Q$ be a Hamiltonian directed path in $T_2-\{a,b,u_1\}$, then $v\, v_1\, b\, a\, u_1\, Q\, v_2\equiv P$. For $(x_{i-1},x_i)\in E(P)$ or $i=1$, if $T_2-u_1\supseteq Q\equiv x_{i+3}...x_n$, then $v_1...v_{i-1}\, v\, u_1\, v_i \, Q\equiv P$. Otherwise, $T_2-u_1\in \mathcal{T}_{3,5,7}$ and $x_{i+3}...x_n$ is antidirected. If $b_2=3$, let $a\in T_2-u_1$, then $T_2-\{a,u_1\}\supseteq Q_1$ such that $v_1...v_{i-1}\,v\,a\, u_1\,v_i\, Q_1\equiv P$. Otherwise, $b_2=2$. If $|T_1|<|T_2|$, then, by considering $\widetilde{P}$ or $\overline{\widetilde{P}}$, we recover one of the previous cases. If $|T_1|=|T_2|$, then the problem is solved unless $T_1$ is isomorphic to $\overline{T_2}$ and $x_1...x_i\equiv x_n...x_{n-i+1}$. Then $\exists\, a\in T_1$ such that $d^+_{T_1}(a)=0$, $T_1-a\in \mathcal{T}_{3,5,7}$ and $N^-(a)\cap T_2=\{u_1\}$. Let $bcd$ be a directed path in $T_1-a$ and let $b'\in T_2-u_1$. $(T_1-\{a,b,c,d\})+v\supseteq Q'$ and $T_2-\{b',u_1\}\supseteq Q''$ such that $Q'\,u_1\,b\,c\,d\,b'\,a\, Q''\equiv P$.\\
From now on, we may suppose that any copy $y_1...y_i$ of $x_1...x_i$ in $T_1$ satisfies the condition that $N^+(y_i)\cap T_2=\phi$. Thus, by simple lemma, we can deduce that $x_i...x_1$ is directed or is of type $P^+(1,2)$ with $T_1-y_i$ is a cyclic triangle. For $i=1$, if $x_{i+2}..x_n$ is not directed, then choose $a\in T_2$ such that $(T_2-a)+v\supseteq Q$ with $a\, v_1\, Q\equiv P$. Otherwise, by considering $\overline{\widetilde{P}}$, the problem is solved since $b_1\geq 2$. For $i\geq2$, let $a\in T_2\cap N^-(v_{i-1})$, then $(T_2-a)+v\supseteq Q$ such that $v_1...v_{i-1}\, a\, v_i\, Q\equiv P$ unless $x_{i+2}...x_n$ is directed or $\in P^+(1,2)$ with $T_2-a$ is a cyclic triangle. If $\exists \, b\in (T_2-a)\cap N^+(a)$, then, if $T_2-a=bcd$ and $x_{i+2}...x_n=P^+(1,2)$, $v_1...v_{i-1}\,a\,b\,d\,v_i\,c\, v\equiv P$ or $x_{i+2}...x_n$ is directed and $v_1...v_{i-1}\,a\,b\,v\,Q\,v_i\equiv P$, where $Q$ is a Hamiltonian directed path in $T_2-\{a,b\}$. If $(T_2-a)\cap N^+(a)=\phi$, since $d^-_{T_2}(v_{i-1})=1$, then $d^+_{T_1}(v_{i-1})=0$ and $i=2$. Thus, by considering $\widetilde{P}=y_1...y_n$ or $\overline{\widetilde{P}}=y_1...y_n$, the problem is solved since $\widetilde{P}$ or $\overline{\widetilde{P}}$ contains the arcs $(y_1,y_2)$ and $(y_2,y_3)$.
\end{case}
\end{proof}
\end{theorem}
\begin{acknowledgment}
I would like to thank Professor Amin El Sahili and Doctor Maydoun Mortada for their following up during the preparation of this paper.
\end{acknowledgment}
\DeclareFieldFormat[article]{pages}{#1}
\renewbibmacro{in:}{%
  \ifentrytype{article}
	{}
	{\bibstring{in}%
	\printunit{\intitlepunct}
	}
}
\DeclareFieldFormat[book]{title}{\mkbibquote{\mkbibemph{#1}}}
\printbibliography[sorting=nty]
\end{document}